\documentclass{amsart}

\usepackage{amssymb}
\usepackage{amsmath}
\usepackage{amsthm}
\usepackage{graphics}

\allowdisplaybreaks

\numberwithin{equation}{section}

\theoremstyle{plain}
\newtheorem{thm}{Theorem}[section]
\newtheorem{prop}[thm]{Proposition}

\newtheorem{cor}[thm]{Corollary}

\theoremstyle{definition}

\newtheorem{rem}[thm]{Remark}

\newcommand{\R}{\mathbb{R}}

\newcommand{\Z}{\mathbb{Z}}

\newcommand{\calF}{\mathcal{F}}

\newcommand{\calS}{\mathcal{S}}

\begin{document}

\title[H\"ormander multiplier theorem and modulation spaces]
{On the H\"ormander multiplier theorem and modulation spaces}
\author{Naohito Tomita}

\address{Naohito Tomita \\
Department of Mathematics \\
Graduate School of Science \\
Osaka University \\
Toyonaka, Osaka 560-0043, Japan}
\email{tomita@math.sci.osaka-u.ac.jp}


\begin{abstract}
It is known that the Sobolev space $L^2_s$ with $s>n/2$
appeared in the H\"ormander multiplier theorem
can be replaced by the Besov space $B^{2,1}_{n/2}$.
On the other hand, the Besov space $B_{n/2}^{2,1}$
is continuously embedded in the modulation space $M^{2,1}_0$.
In this paper, we consider the problem
whether we can replace $B_{n/2}^{2,1}$ by $M^{2,1}_0$.
\end{abstract}
\maketitle

\section{Introduction}\label{section1}
Sj\"ostrand \cite{Sjostrand} proved the $L^2$-boundedness
of pseudo-differential operators with symbols
in the modulation space $M^{\infty,1}(\R^{2n})$
which contains the H\"ormander class $S_{0,0}^0$.
Since then, modulation spaces have been recognized
as a useful tool for pseudo-differential operators.
See B\'enyi-Gr\"ochenig-Okoudjou-Rogers \cite{BGOR},
Cordero-Nicola-Rodino \cite{CNR},
Gr\"ochenig-Heil \cite{Grochenig-Heil} and Toft \cite{Toft}
for further developments.
The purpose of this paper is to apply
modulation spaces to (singular) Fourier multipliers.
\par
We recall some known results on the boundedness of
Fourier multipliers on $L^p(\R^n)$.
The Mihlin multiplier theorem says that
if $m \in C^{[n/2]+1}(\R^{n}\setminus\{0\})$
satisfies
\begin{equation}\label{(1.1)}
|\partial^{\alpha}m(\xi)| \le C_{\alpha}|\xi|^{-|\alpha|}
\quad \text{for all $|\alpha| \le [n/2]+1$}
\end{equation}
then $m(D)$ is bounded on $L^p(\R^n)$ for all $1<p<\infty$
(see \cite[Corollary 8.11]{Duo}),
where $[n/2]$ stands for the largest integer $\le n/2$.
Let $\psi \in \calS(\R^n)$ be such that
$\psi \ge c>0$ on $\{2^{-1/2} \le |\xi| \le 2^{1/2}\}$
and $\mathrm{supp}\, \psi \subset \{2^{-1} \le |\xi| \le 2\}$.
For $m \in \calS'(\R^n)$, we set
\begin{equation}\label{(1.2)}
m_j(\xi)=\psi(\xi)\, m(2^{j}\xi).
\end{equation}
The H\"ormander multiplier theorem \cite{Hormander} states that
if $m \in \calS'(\R^n)$ satisfies
\begin{equation}\label{(1.3)}
\sup_{j \in \Z}\|m_j\|_{L^2_s}<\infty
\quad \text{with} \quad s>n/2
\end{equation}
then $m(D)$ is bounded on $L^p(\R^n)$ for all $1<p<\infty$
(see also \cite[Theorem 8.10]{Duo}),
where $L^2_s(\R^n)$ is the Sobolev space.
We note that \eqref{(1.3)} is weaker than \eqref{(1.1)}.
By using the Besov space $B^{2,1}_{n/2}(\R^n)$
instead of the Sobolev space $L^2_s(\R^n)$ in \eqref{(1.3)},
Seeger \cite{Seeger-2} proved that if $m \in \calS'(\R^n)$ satisfies
\begin{equation}\label{(1.4)}
\sup_{j \in \Z}\|m_j\|_{B^{2,1}_{n/2}}<\infty
\end{equation}
then $m(D)$ is bounded from the Hardy space $H^1(\R^n)$
to the Lorentz space $L^{1,2}(\R^n)$
(see \cite{Seeger-1,Seeger-2} for the definition of $L^{1,2}$).
Then, by interpolation and duality,
\eqref{(1.4)} implies the boundedness of $m(D)$
on $L^p(\R^n)$ for all $1<p<\infty$.
It should be pointed out that
the $L^p$-boundedness of $m(D)$ satisfying \eqref{(1.4)}
follows from a slight modification
of Stein's approach in \cite[Chapter 4, Section 3]{Stein-2}.
Since
\[
L^2_s(\R^n)=B^{2,2}_{s}(\R^n)
\hookrightarrow B^{2,1}_{n/2}(\R^n)
\qquad \text{if $s>n/2$},
\]
we see that \eqref{(1.4)} is weaker than \eqref{(1.3)}.
\par
It is known that the Besov space $B^{2,1}_{n/2}(\R^n)$
is continuously embedded in the modulation space $M^{2,1}_0(\R^n)$,
and this embedding yields the problem
\begin{center}
$\lq\lq$Can we replace $B^{2,1}_{n/2}(\R^n)$ in \eqref{(1.4)}
by $M^{2,1}_0(\R^n)$?".
\end{center}
At least, we have
\begin{thm}\label{1.1}
Let $s>0$.
If $m \in \calS'(\R^n)$ satisfies
\begin{equation}\label{(1.5)}
\sup_{j \in \Z}\|m_j\|_{M^{2,1}_{s}}<\infty,
\end{equation}
then $m(D)$ is bounded on the Hardy space $H^1(\R^n)$,
where $m_j$ is defined by \eqref{(1.2)}.
\end{thm}
We note that,
if $m$ satisfies \eqref{(1.5)} with $s\ge 0$,
then $m(D)$ is bounded on $L^2(\R^n)$
(see the proof of Theorem \ref{1.1}).
Then, by interpolation and duality,
\eqref{(1.5)} with $s>0$
implies the boundedness of $m(D)$ on $L^p(\R^n)$
for all $1<p<\infty$.
Hence, Theorem \ref{1.1} covers
the H\"ormander multiplier theorem, since
\[
L^2_s(\R^n)=M^{2,2}_s(\R^n)
\hookrightarrow M^{2,1}_{s'}(\R^n)
\qquad \text{if $s'<s-n/2$}.
\]
Let us compare \eqref{(1.4)} and \eqref{(1.5)}.
Toft \cite[Theorem 3.1]{Toft} proved the embeddings
\[
B^{2,1}_{n/2}(\R^n) \hookrightarrow M^{2,1}_0(\R^n)
\hookrightarrow B^{2,1}_{0}(\R^n),
\]
and the optimality was proved by
Sugimoto-Tomita \cite[Theorem 1.2]{Sugimoto-Tomita}
(see also \cite{Wang-Huang}).
More precisely,
if $B^{2,1}_s(\R^n) \hookrightarrow M^{2,1}_0(\R^n)$ then $s \ge n/2$,
and if $M^{2,1}_0(\R^n) \hookrightarrow B^{2,1}_s(\R^n)$
then $s \le 0$. 
Then, since $\|f\|_{B^{2,1}_s}\asymp\|(I-\Delta)^{s/2}f\|_{B^{2,1}_0}$
and $\|f\|_{M^{2,1}_s}\asymp\|(I-\Delta)^{s/2}f\|_{M^{2,1}_0}$,
we see that
$B^{2,1}_{n/2}(\R^n) \hookrightarrow M^{2,1}_{s}(\R^n)$
if and only if $s \le 0$,
and $M^{2,1}_{s}(\R^n) \hookrightarrow B^{2,1}_{n/2}(\R^n)$
if and only if $s \ge n/2$.
Therefore, $B^{2,1}_{n/2}(\R^n)$ and $M^{2,1}_s(\R^n)$
have no inclusion relation with each others if $0<s<n/2$:
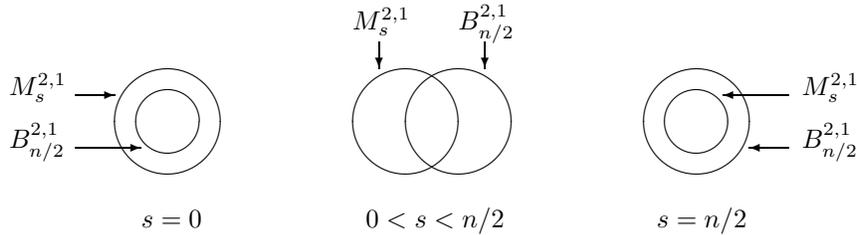
\begin{center}
\begin{picture}(300,100)
\put(40,10){$s=0$}
\put(15,60){\vector(1,0){15}}
\put(-10,60){$M^{2,1}_s$}
\put(15,40){\vector(1,0){25}}
\put(-10,40){$B^{2,1}_{n/2}$}
\put(50,50){\circle{40}}
\put(50,50){\circle{25}}
\put(125,10){$0<s<n/2$}
\put(130,80){\vector(0,-1){10}}
\put(120,85){$M_s^{2,1}$}
\put(140,50){\circle{40}}
\put(170,80){\vector(0,-1){10}}
\put(160,85){$B_{n/2}^{2,1}$}
\put(160,50){\circle{40}}
\put(235,10){$s=n/2$}
\put(285,60){\vector(-1,0){25}}
\put(290,60){$M^{2,1}_s$}
\put(285,40){\vector(-1,0){15}}
\put(290,40){$B^{2,1}_{n/2}$}
\put(250,50){\circle{40}}
\put(250,50){\circle{25}}
\end{picture}
\end{center}
\par
We also mention the relation between
Theorem \ref{1.1} and Baernstein-Sawyer \cite{B-S}
(see also Carbery \cite{Carbery} and Seeger \cite{Seeger}
for some related results).
Since
\begin{equation}\label{(1.6)}
C^{-1}\|\widehat{m_j}\|_{K^{1,1}_s}
\le \|m_j\|_{M^{2,1}_s}
\le C\|\widehat{m_j}\|_{K^{1,1}_s}
\end{equation}
(see Section \ref{section3}),
where $K^{1,1}_s$ is the Herz space,
we have by Theorem \ref{1.1}
\begin{cor}\label{1.2}
Let $s>0$.
If $m \in \calS'(\R^n)$ satisfies
\[
\sup_{j \in \Z}\|\widehat{m_j}\|_{K^{1,1}_s}<\infty,
\]
then $m(D)$ is bounded on the Hardy space $H^1(\R^n)$,
where $m_j$ is defined by \eqref{(1.2)}.
\end{cor}
We remark that Corollary \ref{1.2}
is a special case of \cite[Theorem 3b]{B-S}.
As another corollary of Theorem \ref{1.1},
we have by the norm equivalence
\begin{equation}\label{(1.7)}
C_p^{-1}\|m_j\|_{M^{p,1}_s}
\le \|m_j\|_{M^{2,1}_s}
\le C_p\|m_j\|_{M^{p,1}_s}
\end{equation}
\begin{cor}\label{1.3}
Let $1 \le p \le \infty$ and $s>0$.
If $m \in \calS'(\R^n)$ satisfies
\[
\sup_{j \in \Z}\|m_j\|_{M^{p,1}_s}<\infty,
\]
then $m(D)$ is bounded on the Hardy space $H^1(\R^n)$,
where $m_j$ is defined by \eqref{(1.2)}.
\end{cor}
However, in the critical case $s=0$,
we have the following negative answer:
\begin{prop}\label{1.4}
Let $1<p<\infty$ and $p \neq 2$.
Then there exists a Fourier multiplier $m \in \calS'(\R^n)$ such that
$\sup_{j \in \Z}\|m_j\|_{M^{2,1}_{0}}<\infty$,
but $m(D)$ is not bounded on $L^p(\R^n)$.
\end{prop}
The proofs of Theorem \ref{1.1}, \eqref{(1.6)},
\eqref{(1.7)} and Proposition \ref{1.4} will be given
in Section \ref{section3}.

\section{Preliminaries}\label{section2}
Let $\calS(\R^n)$ and $\calS'(\R^n)$ be the Schwartz spaces of
all rapidly decreasing smooth functions
and tempered distributions,
respectively.
We define the Fourier transform $\calF f$
and the inverse Fourier transform $\calF^{-1}f$
of $f \in \calS(\R^n)$ by
\[
\calF f(\xi)
=\widehat{f}(\xi)
=\int_{\R^n}e^{-i\xi \cdot x}\, f(x)\, dx
\quad \text{and} \quad
\calF^{-1}f(x)
=\frac{1}{(2\pi)^n}
\int_{\R^n}e^{ix\cdot \xi}\, f(\xi)\, d\xi.
\]
For $m \in \calS'(\R^n)$,
the Fourier multiplier operator $m(D)$ is defined by
$m(D)f=\calF^{-1}[m\, \widehat{f}]$ for $f \in \calS(\R^n)$.
The notation $A \asymp B$ stands for $C^{-1}A \le B \le CA$
for some positive constant $C$ independent of $A$ and $B$.
\par
We introduce Besov and modulation spaces,
and suppose that $1 \le p,q \le \infty$ and $s \in \R$.
Let $\psi \in \calS(\R^n)$ be such that
$\psi \ge c>0$ on $\{2^{-1/2}\le |\xi| \le 2^{1/2}\}$,
\begin{equation}\label{(2.1)}
\mathrm{supp}\, \psi \subset \{1/2 \le |\xi| \le 2\}
\quad \text{and} \quad
\sum_{j \in \Z}\psi(2^{-j}\xi)=1
\quad \text{for all $\xi \neq 0$.}
\end{equation}
We set
\begin{equation}\label{(2.2)}
\psi_0(\xi)=1-\sum_{j=1}^{\infty}\psi(2^{-j}\xi)
\quad \text{and} \quad
\psi_j(\xi)=\psi(2^{-j}\xi)
\quad \text{if $j \ge 1$}.
\end{equation}
Then the Besov space $B^{p,q}_s(\R^n)$ consists of
all $f \in \calS'(\R^n)$ such that
\[
\|f\|_{B^{p,q}_s}
=\left(\sum_{j=0}^{\infty}2^{jsq}
\|\psi_j(D)f\|_{L^p}^q\right)^{1/q}<\infty
\]
(with obvious modification in the case $q=\infty$).
We refer to Triebel \cite{Triebel-Book}
and the references therein for more details on Besov spaces.
Let $\varphi \in \calS(\R^n)$ be such that
\begin{equation}\label{(2.3)}
\mathrm{supp}\, \varphi \subset [-1,1]^n
\quad \text{and} \quad
\sum_{k \in \Z^n}\varphi(\xi-k)=1
\quad
\text{for all $\xi \in \R^n$}.
\end{equation}
Then the modulation space $M^{p,q}_s(\R^n)$
consists of all $f \in \calS'(\R^n)$ such that
\[
\|f\|_{M^{p,q}_s}
=\left( \sum_{k \in \Z^n}(1+|k|)^{sq}
\|\varphi(D-k)f\|_{L^p}^q\right)^{1/q}<\infty
\]
(with obvious modification in the case $q=\infty$).
We remark that
\begin{equation}\label{(2.4)}
\|f\|_{M^{p,q}_s}\asymp
\left\{\int_{\R^n}\left(\int_{\R^n}
|V_gf(x,\xi)|^p\, dx\right)^{q/p}(1+|\xi|^2)^{sq/2}\, d\xi\right\}^{1/q},
\end{equation}
where $V_gf$ is the short-time Fourier transform
of $f \in \calS'(\R^n)$
with respect to $g \in \calS(\R^n)\setminus\{0\}$
defined by
\[
V_gf(x,\xi)
=\int_{\R^n}f(t)\, \overline{g(t-x)}\, e^{-i\xi\cdot t}\, dt
\quad \text{for $x,\xi \in \R^n$}
\]
(see, for example, \cite{Triebel}).
The definition of $M^{p,q}(\R^n)$ is independent
of the choice of the window function
$g \in \calS(\R^n)\setminus \{0\}$,
that is,
different window functions
yield equivalent norms
(\cite[Proposition 11.3.2]{Grochenig}).
It is also well known that $M^{2,2}_s(\R^n)=L^2_s(\R^n)$
(\cite[Proposition 11.3.1]{Grochenig}),
where $L^2_s(\R^n)$ is the Sobolev space
defined by the norm
$\|f\|_{L^2_s}=\|(I-\Delta)^{s/2}f\|_{L^2}$
and $(I-\Delta)^{s/2}f=\calF^{-1}[(1+|\xi|^2)^{s/2}\, \widehat{f}]$.
We refer to
Feichtinger \cite{Feichtinger} and Gr\"ochenig \cite{Grochenig}
for more details on modulation spaces
(see also B\'enyi-Grafakos-Gr\"ochenig-Okoudjou \cite{BGGO},
Feichtinger-Narimani \cite{Feichtinger-Narimani}
for Fourier multipliers on modulation spaces).
\par
We next introduce the Hardy and Herz spaces.
Let $\eta \in \calS(\R^n)$ be such that $\int_{\R^n}\eta(x)\, dx=1$.
Then the Hardy space $H^1(\R^n)$ consists of all $f \in L^1(\R^n)$
such that
\[
\|f\|_{H^1}=\int_{\R^n}\sup_{t>0}|\eta_t*f(x)|\, dx<\infty,
\]
where $\eta_t(x)=t^{-n}\eta(t^{-1}x)$.
It is well known that
\begin{equation}\label{(2.5)}
\|f\|_{H^1}
\asymp \|f\|_{L^1}+\sum_{j=1}^n\|R_jf\|_{L^1},
\end{equation}
where $R_j$ is the Riesz transform defined by
\[
R_jf(x)
=\frac{1}{(2\pi)^n}
\int_{\R^n}e^{ix\cdot\xi}
\left(-i\, \frac{\xi_j}{|\xi|}\right)
\widehat{f}(\xi)\, d\xi.
\]
The Herz space $K^{p,q}_s(\R^n)$ consists of
all $f \in L^1_{loc}(\R^n)$ such that
\[
\|f\|_{K^{p,q}_s}
=\left( \sum_{j=0}^{\infty}2^{jsq}
\|\psi_j\, f\|_{L^p}^q \right)^{1/q}<\infty,
\]
where $\{\psi_j\}_{j=0}^{\infty}$ is as in \eqref{(2.2)}.
See Baernstein-Sawyer \cite{B-S} and Stein \cite{Stein}
for more details on Hardy and Herz spaces.

\section{Proof}\label{section3}
Before proving Theorem \ref{1.1},
we note that $M^{2,1}_0 \hookrightarrow \calF L^1 $.
In fact,
by Schwarz's inequality and Plancherel's theorem,
\begin{equation}\label{(1.1.0)}
\begin{split}
\|\widehat{f}\|_{L^1}
&\le \sum_{k \in \Z^n}\|\varphi(\cdot-k)\, \widehat{f}\|_{L^1}
\le C\sum_{k \in \Z^n}\|\varphi(\cdot-k)\, \widehat{f}\|_{L^2}
\\
&=C\sum_{k \in \Z^n}\|\varphi(D-k)\, f\|_{L^2}
=C\|f\|_{M^{2,1}_{0}},
\end{split}
\end{equation}
where $\{\varphi(\cdot-k)\}_{k \in \Z^n}$ is as in \eqref{(2.3)}.
\begin{proof}[Proof of Theorem \ref{1.1}]
Let $\psi$ be as in \eqref{(2.1)} and
$\sup_{j \in \Z}\|m_j\|_{M^{2,1}_{s}}<\infty$,
where $s>0$ and $m_j(\xi)=\psi(\xi)\, m(2^j\xi)$.
Since
\begin{equation}\label{(1.1.1)}
\|m_j\|_{L^{\infty}}
\le C\|\widehat{m_j}\|_{L^1}
\le C\|m_j\|_{M^{2,1}_{0}}
\le C\|m_j\|_{M^{2,1}_{s}}
\end{equation}
and $\psi(\xi)\ge c>0$ on $\{2^{-1/2} \le |\xi| \le 2^{1/2}\}$,
we see that $m \in L^{\infty}$.
This implies that $m(D)$ is bounded on $L^2$.
Then, by the Calder\'on-Zygmund theory
(see, for example,
\cite[Corollary 6.3]{Duo}, \cite[Chapter 3, Theorem 3]{Stein}),
if $K=\calF^{-1}m \in L^1_{loc}(\R^n\setminus\{0\})$ and
\begin{equation}\label{(1.1.2)}
\sup_{y \neq 0}\int_{|x|>2|y|}|K(x-y)-K(x)|\, dx<\infty,
\end{equation}
then $m(D)$ is bounded from $H^1$ to $L^1$.
\par
We only cinsider \eqref{(1.1.2)}
(see Remark \ref{3.1} for the proof of
$K \in L^1_{loc}(\R^n\setminus\{0\})$).
By \eqref{(2.1)}, we have
\[
m(\xi)=\sum_{j \in \Z}\psi(2^{-j}\xi)\, m(\xi)
=\sum_{j \in \Z}m_j(2^{-j}\xi),
\]
and consequently
$K(x)=\sum_{j \in \Z}2^{jn}K_j(2^{j}x)$,
where $K_j=\calF^{-1}m_j$.
Then,
\[
\int_{|x|>2|y|}|K(x-y)-K(x)|\, dx
\le \sum_{j \in \Z}\int_{|x|>2|y|}
|2^{jn}K_j(2^j(x-y))-2^{jn}K_j(2^jx)|\, dx.
\]
Note that $\mathrm{supp}\, m_j \subset \{2^{-1}\le |\xi| \le 2\}$
for all $j \in \Z$.
Since
\begin{equation}\label{(1.1.3)}
\|K_j\|_{L^1}\le C\|m_j\|_{M^{2,1}_0}
\quad \text{and} \quad
\|\nabla K_j\|_{L^1}
\le C\|K_j\|_{L^1} \le C\|m_j\|_{M^{2,1}_0}
\end{equation}
(see \eqref{(1.1.0)} for the left hand inequality,
and \cite[Theorem 1.4.1 (ii)]{Triebel-Book} for the right hand one),
we have by Taylor's formula
\begin{equation}\label{(1.1.4)}
\begin{split}
&\int_{|x|>2|y|}
|2^{jn}K_j(2^j(x-y))-2^{jn}K_j(2^jx)|\, dx
\\
&=\int_{|x|>2|y|}
\left|2^{jn}\sum_{\ell=1}^n (-2^jy_{\ell})
\int_0^1 (\partial_{\ell}K_j)(2^j(x-ty))\, dt \right| dx
\\
&\le C2^j|y|\|\nabla K_j\|_{L^1}
\le C2^j|y|\|m_j\|_{M^{2,1}_s},
\end{split}
\end{equation}
where $y=(y_1,\dots,y_n)$.
On the other hand,
\begin{equation}\label{(1.1.5)}
\int_{|x|>R}|K_j(x)|\, dx \le CR^{-s}\|m_j\|_{M^{2,1}_s}
\quad \text{for all $j \in Z$ and $R>0$}.
\end{equation}
In fact,
since
$\mathrm{supp}\, \varphi(\cdot-k) \subset k+[-1,1]^n
\subset \{|x-k| \le \sqrt{n}\}$ (see \eqref{(2.3)}),
we have
\begin{align*}
&\int_{|x|>R}|K_j(x)|\, dx
\le \sum_{k \in \Z^n}
\int_{|x|>R}|\varphi(-x-k)\, K_j(x)|\, dx
\\
&\le \sum_{|k|>R/2}
\int_{\R^n}|\varphi(-x-k)\, \calF^{-1}{m_j}(x)|\, dx
\le \sum_{|k|>R/2}
\left(\int_{\R^n}|\varphi(x-k)\, \widehat{m_j}(x)|^2\, dx\right)^{1/2}
\\
&=(2\pi)^{n/2}\sum_{|k|>R/2}
(1+|k|)^{-s}(1+|k|)^{s}\|\varphi(D-k)m_j\|_{L^2}
\le CR^{-s}\|m_j\|_{M_s^{2,1}}
\end{align*}
for all $R>2\sqrt{n}$, and
\[
\int_{|x|>R}|K_j(x)|\, dx
\le (1+R)^{-s}(1+R)^s\|K_j\|_{L^1}
\le CR^{-s}\|m_j\|_{M^{2,1}_s}
\]
for all $0<R \le 2\sqrt{n}$,
where we have used \eqref{(1.1.3)}.
Then, \eqref{(1.1.5)} gives
\begin{equation}\label{(1.1.6)}
\begin{split}
&\int_{|x|>2|y|}
|2^{jn}K_j(2^j(x-y))-2^{jn}K_j(2^jx)|\, dx
\\
&\le 2\int_{|x|>|y|}2^{jn}|K_j(2^jx)|\, dx
=2\int_{|x|>2^j|y|}|K_j(x)|\, dx \le C(2^j|y|)^{-s}\|m_j\|_{M^{2,1}_s}
\end{split}
\end{equation}
for all $y \neq 0$.
Hence, it follows from \eqref{(1.1.4)} and \eqref{(1.1.6)} that
\begin{align*}
&\sum_{j \in \Z}\int_{|x|>2|y|}
|2^{jn}K_j(2^j(x-y))-2^{jn}K_j(2^jx)|\, dx
\\
&=\left(\sum_{2^j|y| \le 1}+\sum_{2^j|y|>1}\right)
|2^{jn}K_j(2^j(x-y))-2^{jn}K_j(2^jx)|\, dx
\\
&\le C\sum_{2^j|y| \le 1}2^j|y|\|m_j\|_{M^{2,1}_s}
+C\sum_{2^j|y|>1}(2^j|y|)^{-s}\|m_j\|_{M^{2,1}_s}
\le C\|m_j\|_{M^{2,1}_s}
\end{align*}
for all $y \neq 0$.
Therefore, $m(D)$ is bounded from $H^1$ to $L^1$.
This implies the boundedness of $m(D)$ on $H^1$.
In fact, by \eqref{(2.5)},
\begin{align*}
\|m(D)f\|_{H^1}
&\le C\left(\|m(D)f\|_{L^1}+
\sum_{j=1}^{n}\|R_j(m(D)f)\|_{L^1}\right)
\\
&\le C\left(\|f\|_{H^1}+
\sum_{j=1}^{n}\|R_jf\|_{H^1}\right)
\le C\|f\|_{H^1},
\end{align*}
where we have used the fact that $R_j$ is bounded on $H^1$.
The proof is complete.
\end{proof}
\begin{rem}\label{3.1}
It is not difficult to prove that,
if $m$ satisfies \eqref{(1.5)} with $s>0$,
then $K=\calF^{-1}m \in L^1_{loc}(\R^n\setminus\{0\})$.
Let $m_j$ and $K_j$ be as in the proof of Theorem \ref{1.1},
and recall that
$K(x)=\sum_{j \in \Z}2^{jn}K_j(2^j x)$.
Since $\mathrm{supp}\, m_j \subset \{1/2 \le |\xi| \le 2\}$
for all $j \in \Z$,
it follows from \eqref{(1.1.1)} that
$\|K_j\|_{L^{\infty}}
=\|\calF^{-1}m_j\|_{L^{\infty}}
\le C\|m_j\|_{L^{\infty}} \le C\|m_j\|_{M^{2,1}_s}$
for all $j \in \Z$.
Hence, for any $0<R_1<R_2<\infty$,
\begin{align*}
&\sum_{j=-\infty}^{0}
\int_{R_1 \le |x| \le R_2}
|2^{jn}K_j(2^j x)|\, dx
\le C_{R_1,R_2}\sum_{j=-\infty}^0
2^{jn}\|K_j\|_{L^{\infty}}
\\
&\le C_{R_1,R_2}\sum_{j=-\infty}^0
2^{jn}\|m_j\|_{M^{2,1}_s}
\le C_{R_1,R_2}
\sup_{j \le 0}\|m_j\|_{M^{2,1}_s}.
\end{align*}
On the other hand,
by \eqref{(1.1.5)},
\begin{align*}
&\sum_{j=1}^{\infty}
\int_{R_1 \le |x| \le R_2}
|2^{jn}K_j(2^j x)|\, dx
\le \sum_{j=1}^{\infty}
\int_{|x| \ge 2^j R_1}
|K_j(x)|\, dx
\\
&\le \sum_{j=1}^{\infty}
C(2^jR_1)^{-s}\|m_j\|_{M^{2,1}_s}
\le C_{R_1}\sup_{j \ge 1}\|m_j\|_{M^{2,1}_s}.
\end{align*}
Therefore, we see that
\[
\int_{R_1 \le |x| \le R_2}
|K(x)| \, dx
\le \sum_{j \in \Z}
\int_{R_1 \le |x| \le R_2}
|2^{jn}K(2^j x)| \, dx
\le C_{R_1,R_2}\sup_{j \in \Z}\|m_j\|_{M^{2,1}_s},
\]
that is, $K \in L^1_{loc}(\R^n\setminus\{0\})$.
\end{rem}
To prove \eqref{(1.6)} and \eqref{(1.7)},
we use the following fact \cite[Remark 4.2]{RSTT}
(see also \cite[Lemma 1]{Okoudjou} for the case $s=0$),
and give the proof for reader's convenience.
\begin{prop}\label{3.2}
Let $1 \le p,q \le \infty$, $s \in \R$,
and let $\Omega$ be a compact subset of $\R^n$.
Then there exists a constant $C_{\Omega}>0$ such that
\[
C_{\Omega}^{-1}\|(I-\Delta)^{s/2}f\|_{\calF L^q}
\le \|f\|_{M^{p,q}_s}
\le C_{\Omega}\|(I-\Delta)^{s/2}f\|_{\calF L^q}
\]
for all $f \in \calS'(\R^n)$ with $\mathrm{supp}\, f \subset \Omega$,
where $\|f\|_{\calF L^q}=\|\widehat{f}\|_{L^q}$.
\end{prop}
\begin{proof}
Our proof is based on one of \cite[Lemma 1]{Okoudjou}.
Let $\Omega \subset \{|x| \le R\}$,
and let $f \in \calS'(\R^n)$ and
$g \in \calS(\R^n)\setminus\{0\}$ be such that
$\mathrm{supp}\, f \subset \{|x| \le R\}$,
$\mathrm{supp}\, g \subset \{|x| \le 4R\}$
and $g=1$ on $\{|x| \le 2R\}$.
Then
$\mathrm{supp}\, V_gf(\cdot,\xi)
\subset \{x: |x| \le 5R\}$
for all $\xi \in \R^n$.
By Plancherel's theorem,
\[
\left|V_gf(x,\xi)\right|
=\frac{1}{(2\pi)^n}
\left|\int_{\R^n}\widehat{f}(t)\,
\overline{\widehat{g}(t-\xi)}\, e^{ix\cdot t}\, dt\right|
\le \frac{1}{(2\pi)^n}
\int_{\R^n}\left|\widehat{f}(t)\,
\widehat{g}(t-\xi)\right| dt
\]
for all $x \in \R^n$.
Hence, by \eqref{(2.4)},
\begin{align*}
&\|f\|_{M^{p,q}_s}
\le C\left[ \int_{\R^n} \left\{
\int_{|x| \le 5R}
\left( \int_{\R^n}\left|\widehat{f}(t)\,
\widehat{g}(t-\xi)\right| dt \right)^p dx
\right\}^{q/p} (1+|\xi|^2)^{sq/2}\, d\xi \right]^{1/q}
\\
&\le CR^{n/p}\left\{\int_{\R^n}
\left( \int_{\R^n}(1+|t|^2)^{s/2}\,|\widehat{f}(t)|\,
(1+|t-\xi|^2)^{|s|/2}\, |\widehat{g}(t-\xi)|\,
dt \right)^q d\xi \right\}^{1/q}
\\
&\le C_R\|(I-\Delta)^{|s|/2}g\|_{\calF L^1}
\|(I-\Delta)^{s/2}f\|_{\calF L^q},
\end{align*}
where $C>0$ is independent of $f$.
\par
On the other hand,
since $g=1$ on $\{|x| \le 2R\}$,
we see that
\[
\widehat{f}(\xi)
=\int_{\R^n}f(t)\, e^{-i\xi\cdot t}\, dt
=\int_{|t| \le R}f(t)\, \overline{g(t-x)}\, e^{-i\xi\cdot t}\, dt
=V_gf(x,\xi)
\]
for all $|x| \le R$.
This gives
\[
|\widehat{f}(\xi)| \le CR^{-n/p}
\left(\int_{|x| \le R}|V_gf(x,\xi)|^p\, dx \right)^{1/p}
\le CR^{-n/p}\|V_gf(\cdot, \xi)\|_{L^p}
\]
for all $\xi \in \R^n$.
Therefore,
\begin{align*}
&\|(I-\Delta)^{s/2}f\|_{\calF L^q}
=\left( \int_{\R^n}|(1+|\xi|^2)^{s/2}
\widehat{f}(\xi)|^q\, d\xi \right)^{1/q}
\\
&\le C_R\left\{ \int_{\R^n}\left((1+|\xi|^2)^{s/2}
\|V_gf(\cdot, \xi)\|_{L^p}\right)^q d\xi \right\}^{1/q}
\le C_R\|f\|_{M^{p,q}_s},
\end{align*}
where $C>0$ is independent of $f$.
The proof is complete.
\end{proof}
We are now ready to prove \eqref{(1.6)}, \eqref{(1.7)}
and Proposition \ref{1.4}.
\begin{proof}[Proofs of \eqref{(1.6)} and \eqref{(1.7)}]
Let $m_j$ be defined by \eqref{(1.2)}.
Note that $\mathrm{supp}\, m_j \subset \{2^{-1}\le |\xi| \le 2\}$
for all $j \in \Z$.
Then, by Proposition \ref{3.2},
\begin{equation}\label{(1.6.1)}
C^{-1}\|(I-\Delta)^{s/2}m_j\|_{\calF L^1}
\le \|m_j\|_{M^{2,1}_s}
\le C\|(I-\Delta)^{s/2}m_j\|_{\calF L^1}
\end{equation}
for all $j \in \Z$.
On the other hand,
since $(1+|x|^2)^{1/2} \asymp 2^{\ell}$
on $\mathrm{supp}\, \psi_{\ell}$ for all $\ell \ge 0$
(see \eqref{(2.1)} and \eqref{(2.2)}),
we have
\begin{equation}\label{(1.6.2)}
\begin{split}
\|(I-\Delta)^{s/2}m_j\|_{\calF L^1}
&\asymp \sum_{\ell=0}^{\infty}
\int_{\R^n}\left|\psi_{\ell}(x)\, (1+|x|^2)^{s/2}\,
\widehat{m_j}(x)\right| dx
\\
&\asymp \sum_{\ell=0}^{\infty}2^{\ell s}
\int_{\R^n}\left|\psi_{\ell}(x)\, \widehat{m_j}(x)\right| dx
=\|\widehat{m_j}\|_{K^{1,1}_s}.
\end{split}
\end{equation}
Hence, combining \eqref{(1.6.1)} and \eqref{(1.6.2)},
we have \eqref{(1.6)}.
\par
Let $1 \le p \le \infty$.
Then, by Proposition \ref{3.2},
\begin{equation}\label{(1.7.1)}
C_p^{-1}\|(I-\Delta)^{s/2}m_j\|_{\calF L^1}
\le \|m_j\|_{M^{p,1}_s}
\le C_p\|(I-\Delta)^{s/2}m_j\|_{\calF L^1}
\end{equation}
for all $j \in \Z$.
Therefore, combining \eqref{(1.6.1)} and \eqref{(1.7.1)},
we have \eqref{(1.7)}.
\end{proof}
\begin{proof}[Proof of Proposition \ref{1.4}]
The following counterexample
(Triebel \cite[Proposition 2.6.4]{Triebel-Book}) is known:
\begin{equation}\label{(1.4.1)}
\begin{cases}
\text{$m(D)$ is bounded on $B^{p,q}_s(\R^n)$
for all $1 \le p,q \le \infty$ and $s \in \R$}
\\
\text{$m(D)$ is not bounded on $L^{p}(\R^n)$ for any $p \neq 2$}
\end{cases}
\end{equation}
(see also Littman-McCarthy-Rivi\`ere \cite{L-Mc-R}
and Stein-Zygmund \cite{Stein-Zygmund}).
Let $m$ be as in \eqref{(1.4.1)},
and we prove that $m$ satisfies
$\sup_{j \in \Z}\|m_j\|_{M^{2,1}_0}<\infty$,
where $m_j$ is defined by \eqref{(1.2)}.
We remark that $m(D)$ is bounded on $B^{1,q}_s(\R^n)$
for some $1 \le q \le \infty$ and $s \in R$
if and only if $\calF^{-1}m \in B^{1,\infty}_0(\R^n)$
(see \cite[Theorem 2.6.3]{Triebel-Book}).
Then, the boundedness of $m(D)$ on $B^{1,q}_s$ implies
$\calF^{-1}m \in B^{1,\infty}_0$.
Hence,
since $B^{1,\infty}_0 \hookrightarrow \Dot{B}^{1,\infty}_0$,
we see that
\begin{equation}\label{(1.4.2)}
\begin{split}
&\sup_{j \in \Z}\|\calF^{-1}m_j\|_{L^1}
=\sup_{j \in \Z}
\|2^{-jn}\calF^{-1}[\psi(2^{-j}\cdot)\, m](2^{-j}\cdot)\|_{L^1}
\\
&=\sup_{j \in \Z}\|\calF^{-1}[\psi(2^{-j}\cdot)\, m]\|_{L^1}
=\sup_{j \in \Z}\|\psi(2^{-j}D)(\calF^{-1}m)\|_{L^1}
\\
&=\|\calF^{-1}m\|_{\Dot{B}^{1,\infty}_0}
\le C\|\calF^{-1}m\|_{B^{1,\infty}_0}<\infty.
\end{split}
\end{equation}
On the other hand,
since $\mathrm{supp}\, m_j \subset \{2^{-1} \le |\xi| \le 2\}$
for all $j \in \Z$, we have by Proposition \ref{3.2}
\begin{equation}\label{(1.4.3)}
C^{-1}\|\calF^{-1}m_j\|_{L^1}
\le \|m_j\|_{M^{2,1}_0} \le C\|\calF^{-1}m_j\|_{L^1}
\quad \text{for all $j \in \Z$.}
\end{equation}
Therefore, combining \eqref{(1.4.2)} and \eqref{(1.4.3)},
we see that $\sup_{j \in \Z}\|m_j\|_{M^{2,1}_0}<\infty$,
but $m(D)$ is not bounded on $L^{p}(\R^n)$ for any $p \neq 2$.
The proof is complete.
\end{proof}

\section*{Acknowledgement}
The author gratefully acknowledges helpful discussions
with Professors Akihiko Miyachi, Michael Ruzhansky,
Mitsuru Sugimoto and Joachim Toft.
He also would like to thank the referees
for their valuable comments and remarks.


\end{document}